\documentclass{amsart}
 
\usepackage[T1]{fontenc}
\usepackage[utf8]{inputenc}
\usepackage[pdftex]{graphicx}
\usepackage{tikz-cd}
\usepackage{blindtext}

\usepackage{hyperref,xcolor}
\usepackage{amsmath,amsfonts,amssymb,amsthm,color,mathtools,enumitem}
\usepackage{epstopdf}
\usepackage{polynom}
\usepackage{babel}
\usepackage{mathrsfs}
\usepackage{chngcntr}
\usepackage{verbatim} 

\hypersetup{
	colorlinks=true,
	linkcolor=blue,
	citecolor=cyan,
	urlcolor=cyan}
	
	\usepackage{hyphenat}

\newtheorem{theorem}{Theorem}[section] 

\newtheorem{cor}[theorem]{Corollary}
\newtheorem{lemma}[theorem]{Lemma}
\newtheorem{example}[theorem]{Example}
\newtheorem{remark}[theorem]{Remark}

\newtheorem*{result*}{Theorem}
\newtheorem*{resultcor*}{Corollary}
\newtheorem{result}{Theorem}

\newtheorem{resultcor}[result]{Corollary}

\numberwithin{equation}{section}

\newcommand{\ee}{\varepsilon}

\newcommand{\Tr}{\text{Tr}}

\newcommand{\ev}{\text{ev}}

\newcommand{\tD}{\tilde{\Delta}}

\newcommand{\tr}{\text{tr}}
\newcommand{\id}{\text{id}}

\newcommand{\II}{\text{II}}

\DeclareMathOperator{\Aff}{Aff}

\newcommand{\MvN}{\sim_{\text{MvN}}}

\newcommand{\dlim}{\underset{\to}{\lim}}

\newcommand{\bN}{\mathbb{N}}
\newcommand{\bR}{\mathbb{R}}

\newcommand{\bT}{\mathbb{T}}

\newcommand{\cM}{\mathcal{M}}

\newcommand{\cZ}{\mathcal{Z}}

\begin{document}
\title[Universal covering groups of unitary groups of W*-algebras]{Universal covering groups of unitary groups of von Neumann algebras}

\author[Pawel Sarkowicz]{Pawel Sarkowicz}
\address[Pawel Sarkowicz]{Department of Pure Mathematics, University of Waterloo,  N2L 3G1, Canada}
\email{\href{mailto:psarkowi@uwaterloo.ca}{psarkowi@uwaterloo.ca}}

  \begin{abstract}
  We give a short and simple proof, utilizing the pre-determinant of P. de la Harpe and G. Skandalis, that the universal covering group of the unitary group of a $\II_1$ von Neumann algebra $\cM$, when equipped with the norm topology, splits algebraically as the direct product of the self-adjoint part of its center and the unitary group $U(\cM)$. Thus, when $\cM$ is a $\II_1$ factor, the universal covering group of $U(\cM)$ is algebraically isomorphic to the direct product $\bR \times U(\cM)$.
  In particular, the question of P. de la Harpe and D. McDuff of whether the universal cover of $U(\cM)$ is a perfect group is answered in the negative.
  \end{abstract}

  \maketitle
  \tableofcontents
  
 \section{Introduction}
  
The universal covering group of a topological group, recoverable as the space of paths modulo null-homotopy when it is sufficiently connected, gives rise to an extension of the group by its fundamental group. 
  The universal covering group of $U(n)$, the unitary group of the $n \times n$ matrices, is known to homeomorphically isomorphic to $SU(n) \rtimes \bR$, where $SU(n)$ is the special unitary group consisting of unitaries with trivial determinant.
Indeed, $U(n)$ itself is the semidirect product $SU(n) \rtimes \bT$ as can see the from the topologically split extension
\begin{equation}\label{eq:U(n) splitting}
\begin{tikzcd}
1 \arrow[r] & SU(n) \arrow[r] & U(n) \arrow[r, "\det"] & \bT \arrow[l, bend right=49] \arrow[r] & 1.
\end{tikzcd}
\end{equation}
Therefore $U(n)$ is covered by $SU(n) \times \bT$, which is further covered by $SU(n) \times \bR$, and this must necessarily be the universal cover, as a topological space, since it is simply connected.
The isomorphism of the universal cover, as a topological group, with the semidirect product $SU(n) \rtimes \bR$ was given in \cite{AguilarSocolovsky00}.
  
  If $\cM$ is an infinite von Neumann factor, then it is known that the universal covering group is in fact just $U(\cM)$ due to the fact that the fundamental group will be trivial -- see \cite{Kuiper65,Breuer70,ArakiSmithSmith71,BruningWillgerodt76,Handelman78,Schroder84,Schroder84real} for results on homotopy groups related to von Neumann algebras.
  Therefore to understand the universal covering group of unitary groups of von Neumann factors, it remains to study them in the $\II_1$ setting.
  We show that the universal covering group of $U(\cM)$, where $\cM$ is any $\II_1$ von Neumann algebra, exhibits similar and even simpler algebraic behaviour to that of the universal covering group of $U(n)$.

\begin{result}
Let $\cM$ be a type $\II_1$ von Neumann algebra. Then the universal covering group of $U(\cM)$ exists and splits algebraically as the direct product
\begin{equation}
\widetilde{U(\cM)} \simeq Z(\cM)_{sa} \times U(\cM).
\end{equation}
In particular, when $\cM$ is a $\II_1$ factor,
\begin{equation}
\widetilde{U(\cM)} \simeq \bR \times U(\cM).
\end{equation}
\end{result}

We however note that this is not a topological product. Indeed, $\pi_1(U(\cM))$ is viewed as a discrete subset. Even if $\bR$ was equipped with its standard topology, $U(\cM)$ is not simply connected in the $\II_1$ setting, so that the product space would not be. 

One observation is that there is a way to define a determinant map on the group of invertibles of a $\II_1$ factor $\cM$ via the Fuglede-Kadison determinant \cite{FugledeKadison52} or the de la Harpe-Skandalis determinant \cite{dlHS84a} (they are related -- see \cite{dlHarpe13} for exposition).
Restricting this map to the unitary group, it is a feature that this restriction gives the trivial map.
In particular, if one were to define the \emph{special unitary group} of $\cM$ to be those unitaries in $\cM$ which have trivial determinant, we would have that $SU(\cM) = U(\cM)$, and so the splitting above can be viewed as the direct product of $\bR$ with the special unitary group of $\cM$. 
In this sense, this theorem above is an algebraic generalization of the fact that the unitary group $U(n)$ has covering group $SU(n) \rtimes \bR$. 

In \cite{dlHarpeMcDuff83}, de la Harpe and McDuff examined the homology groups of certain automorphism groups arising in operator algebras. 
Among other things, they showed that the unitary groups of properly infinite von Neumann algebras are acyclic -- that is, all of their integral homology groups vanish.
The unitary group of a $\II_1$ factor $\cM$ has vanishing first integral homology group (since it is perfect -- see \cite{Broise67,FackdlHarpe80}), but it came into question whether or not the second integral homology group of $U(\cM)$ vanishes. 
If one could show that the universal covering group $\widetilde{U(\cM)}$ was perfect, then this group would be a cover in the algebraic sense \cite{Kervaire70} and there would be potential for non-triviality of the second homology group.
However, as a consequence of the main theorem, we answer the question about the universal covering being perfect in the negative. 

\begin{resultcor}
Let $\cM$ be a $\II_1$ factor. Then the universal covering group of $U(\cM)$ is not a perfect group.
\end{resultcor}
The question of the whether or not the second integral homology group (or all of the higher integral homology groups) of $U(\cM)$ is trivial, where $\cM$ is a $\II_1$ factor, remains open.

  \addtocontents{toc}{\protect\setcounter{tocdepth}{0}}
  \section*{Acknowledgements}
  Thanks to Mehdi Moradi for insightful conversations, as well as to Aaron Tikuisis and Laurent Marcoux for looking over a first draft.

  \addtocontents{toc}{\protect\setcounter{tocdepth}{1}}

\section{Preliminaries}\label{sec:prelim}
\subsection{Notation}
For a group $G$, we denote by $DG$ the derived subgroup of $G$, i.e.,
\begin{equation}
DG := \langle ghg^{-1}h^{-1} \mid g,h \in G \rangle
\end{equation}
is the subgroup generated by commutators.
If $G$ has an underlying topology, we denote by $G^0$ the connected component of the identity, and we write
\begin{equation}
PG= \{\xi: [0,1] \to G \mid \xi \text{ is continuous, } \xi(0) = 1_G \}
\end{equation}
to be the path group. Multiplication is given by the pointwise product and the topology is the compact-open topology, meaning the topology generated by subbase consisting of sets
\begin{equation}
V(K,O) = \{ f \in C([0,1],G) \mid f(K) \subseteq O\},
\end{equation}
where $K \subseteq [0,1]$ is closed and $O \subseteq G$ is open.
For two paths $\xi,\eta \in PG$, we will write $\xi \sim_{nh} \eta$, if $\xi$ and $\eta$ are \emph{null-homotopic}, meaning that the path $\xi^{-1}\eta$ is homotopic to the constant path at our base point (in this case, the identity of $G$). In particular, two null-homotopic paths must have the same endpoint.

For a unital C*-algebra $A$, $A_{sa}$ will denote the self-adjoint elements of $A$, while $U(A)$ will denote the group of unitaries in $A$: that is, the subgroup of invertible elements $u \in A$ satisfying
\begin{equation}
u^*u = uu^* = 1_A.
\end{equation}
We give $U(A)$ the subspace topology induced by the norm of $A$. 
For $n \in \bN$, we let $U_n(A) = U(M_n(A))$ and take
\begin{equation}
U_{\infty}(A) = \dlim U_n(A)
\end{equation}
to be the direct limit with connecting maps $u \mapsto u \oplus 1$.
We equip $U_{\infty}(A)$ with the norm topology.\footnote{One can also equip $U_{\infty}(A)$ with the inductive limit topology. However, even though $U_n(A)$ are topological groups, $U_{\infty}(A)$ is almost never a topological group because multiplication will not be jointly continuous in general. The quotient of the closure of the derived group can be more useful (for example, in terms of classification of morphisms) when one uses the inductive limit topology as opposed to the norm topology -- see \cite{Thomsen95,NielsenThomsen96,CGSTW23}.}
The topological $K_1$ group of a unital C*-algebra $A$ is then the set of connected components of $U_{\infty}(A)$
\begin{equation}
K_1(A) := U_{\infty}(A)/U_{\infty}^0(A).
\end{equation}
For two projections $p,q$ over matrix algebras over $A$, we write $p \MvN q$ if they are Murray-von Neumann equivalent. 
The $K_0$ group of a unital C*-algebra is the Grothendieck group of the Murray-von Neumann semigroup of equivalence classes of projections in matrix algebras over $A$, and there is a canonical isomorphism
\begin{equation}
K_0(A) \simeq \pi_1(U_{\infty}(A))
\end{equation}
which is induced by the map which takes a projection $p \in M_n(A)$ to the loop
\begin{equation}
t \mapsto e^{2\pi i t}p + 1-p
\end{equation}
in $U_n(A) \subseteq U_{\infty}(A)$.
Thus we have
\begin{equation}
K_n(A) = \begin{cases}
\pi_1(U_{\infty}(A)) & n =0; \\
\pi_0(U_{\infty}(A)) & n =1.
\end{cases}
\end{equation}

\subsection{The pre-determinant}
We will shortly recall the definition of the \emph{pre-determinant} map, which will take a piece-wise smooth path to an element of a real Banach space, when we have a \emph{tracial} map.
To extend its definition to any continuous path, we will need a few results.
The following can be found in the proof of \cite[Lemma 3]{dlHS84a}.

\begin{lemma}\label{lem:any path homotopic to piecewise smooth}
Any path $\xi \in PU^0(A)$ is homotopic to a piece-wise smooth.
\end{lemma}
In fact, any path $\xi \in PU^0(A)$ is homotopic to a path of the form
\begin{equation}
t \mapsto \xi\left(\frac{j-1}{k}\right)e^{2\pi i(kt -j + 1)a_j}, t \in \left[\frac{j-1}{k},\frac{j}{k}\right], j=1,\dots,k,
\end{equation}
where $k$ is large enough such that
\begin{equation}
\left\| \xi\left(\frac{j-1}{k}\right)^{-1}\xi\left(\frac{j}{k}\right) - 1\right\| < 1 \text{ for all } j=1,\dots,k
\end{equation}
and
\begin{equation}
a_j = \frac{1}{2\pi i}\log \left(\xi\left(\frac{j-1}{k}\right)^{-1}\xi\left(\frac{j}{k}\right)\right) \text{ for } j=1,\dots,k.
\end{equation}
We make the historical note that some of the ideas of working with paths in this way can be found in \cite{ArakiSmithSmith71}.

If $\tau: A_{sa} \to E$ is a bounded tracial map to a real Banach $E$ (tracial meaning that $\tau(aa^*) = \tau(a^*a)$ for all $a \in A$), we can associate to a piece-wise smooth path $\xi \in PU^0(A)$ the element
\begin{equation}
\tD_\tau(\xi) := \int_0^1\tau\left(\frac{1}{2\pi i}\xi'(t)\xi(t)^{-1}\right)dt \in E,
\end{equation}
which is just the Riemann path-integral in $E$.
We note that this is well-defined since it is easily seen that, when $\xi'(t)$ is defined, $\xi'(t)\xi(t)^{-1}$ is skew-adjoint, so that $\frac{1}{2\pi i}\xi'(t)\xi(t)^{-1}$ is self-adjoint.
This map also extends to $U_n^0(A)$, hence $U_{\infty}^0(A)$, in the canonical way via the unnormalized trace: for $\xi \in PU_n^0(A)$ piece-wise smooth, write
\begin{equation}
\tD_\tau(\xi) := \int_0^1 \tr_n \otimes \tau \left(\frac{1}{2\pi i}\xi'(t)\xi(t)^{-1}\right)dt,
\end{equation}
where $\tr_n \otimes \tau: M_n(A)_{sa} \to E$ is the (unnormalized) extension of the trace, given by $\tr_n\otimes \tau((a_{ij})) = \sum_i \tau(a_{ii})$, whenever $(a_{ij}) \in M_n(A)_{sa}$.
In this case, we have
\begin{equation}
\tD_\tau(\xi \oplus 1_{m-n}) = \tD_\tau(\xi).
\end{equation}
whenever $\xi \in PU_n^0(A)$ and $m > n$.

We state the unitary variant of \cite[Lemme 1]{dlHS84a}, suited to our particular needs.

\begin{lemma}\label{lem:pre-det facts}
Let $\tau: A_{sa} \to E$ be a bounded trace on a unital C*-algebra $A$. The map $\tD$, which takes a piece-wise smooth path $\xi \in PU^0(A)$ (or $\xi \in PU_{\infty}^0(A)$) to $\tD_\tau(\xi)$, has the following properties.
\begin{enumerate}
\item It takes point-wise products to sums: if $\xi_1,\xi_2$ are two piece-wise smooth paths, then
\begin{equation}
\tD_\tau(\xi_1\xi_2) = \tD_\tau(\xi_1) + \tD_\tau(\xi_2),
\end{equation}
where $\xi_1\xi_2$ is the piece-wise smooth path $[\xi_1\xi_2](t) = \xi_1(t)\xi_2(t)$. 
\item If $\|\xi(t) - 1\| < 1$ for all $t \in [0,1]$, then
\begin{equation}
\tD_\tau(\xi) = \frac{1}{2\pi i}\tau\left(\log \xi(1)\right).
\end{equation}
\item $\tD_\tau(\xi)$ depends only on the continuous homotopy class of $\xi$ with fixed endpoints.
\item If $p \in A$ is a projection, then the path $\xi_p \in PU^0(A)$ given by
\begin{equation}
\xi_p(t) = pe^{2\pi i t} + (1-p)
\end{equation}
satisfies $\tD_\tau(\xi_p) = \tau(p)$. Note that $\xi_p(t)$ can also be written as $\xi_p(t) = e^{2\pi i tp}$. 
\end{enumerate}
\end{lemma}

By combining Lemma \ref{lem:any path homotopic to piecewise smooth} and Lemma \ref{lem:pre-det facts}(3), it makes sense to extend $\tD_\tau$ to a map
\begin{equation}
\tD_\tau: PU^0(A) \to E
\end{equation}
by $\tD_\tau(\xi) = \tD_\tau(\xi_0)$ where $\xi_0 \in PU^0(A)$ is any piece-wise smooth path homotopic to $\xi$, and this map is clearly a group homomorphism. 
Lemma \ref{lem:pre-det facts}(4) in fact allows us to determine the pairing map between the $K_0$-group and traces.

We will be interested in $\II_1$ von Neumann algebras $\cM$, which have a normal, faithful, positive, center-valued trace \cite[Theorem V.2.6]{TakesakiI}, which we will write as $\Tr_\cM: \cM \to Z(\cM)$. Moreover, as it is self-adjoint in the sense that $\Tr_\cM(\cM_{sa}) \subseteq Z(\cM)_{sa}$, we will abuse notation and write $\Tr_\cM$ for the restriction to $\cM_{sa}$. We will just write $\tD$ for $\tD_{\Tr_\cM}$.
We would like a generalization of the fact that the pairing between $K$-theory and traces of a $\II_1$ factor $\cM$ gives an isomorphism $K_0(\cM) \simeq \bR$.
We note that, for example in \cite{Handelman78}, it was shown that $\pi_1(U(\cM)) \simeq Z(\cM)_{sa}$ in the $\II_1$ von Neumann algebra setting (actually this was shown to hold more generally). We will want to work with a specific isomorphism.

\begin{lemma}\label{lem:fundamental group is center}
Let $\cM$ be a type $\II_1$ von Neumann algebra. Then
\begin{equation}
\tD|_{\pi_1(U(\cM))}: \pi_1(U(\cM)) \to Z(\cM)_{sa}
\end{equation}
is a group isomorphism. (We are abusing notation by writing $\tD|_{\pi_1(U(\cM))}$, but as $\tD$ is homotopy invariant, it makes sense to apply it to homotopy classes.)
\begin{proof}
As $\cM$ is a $\II_1$ von Neumann algebra, it follows from Corollary 4.5 and Theorem 4.13 of \cite{Rieffel87} that the canonical map
\begin{equation}\label{eq:canonical map is iso}
\pi_1(U(\cM)) \to \pi_1(U_{\infty}(\cM)) = K_0(\cM)
\end{equation}
is an isomorphism. 

Now given some positive contraction $a \in Z(\cM)$, there is a projection $p \in \cM$ such that $\Tr_\cM(p) = a$ by \cite[Theorem 8.4.4(ii)]{KadisonRingroseII}.
As such, this implies that given $n \in \bN$ and any positive element $a \in Z(\cM)$ with $\|a\| \leq n$, we can find a projection in $p \in M_n(\cM)$ such that $\tr_n \otimes \Tr_\cM(p) = a$. 
As every self-adjoint element in $Z(\cM)$ is the difference of positive elements in $Z(\cM)$ and Murray-von Neumann equivalent projections $p,q \in M_n(\cM)$ will satisfy $\Tr_\cM(p) = \Tr_\cM(q)$, we get a well-defined group homomorphism
\begin{equation}
K_0(\Tr_\cM): K_0(\cM) \to Z(\cM)_{sa}
\end{equation}
which is surjective.
To show that this map is injective, we simply note that that $p \MvN q$ if and only if $\Tr_\cM(p) = \Tr_\cM(q)$ by \cite[Theorem 8.4.3(iii)]{KadisonRingroseII}.

Now Lemma \ref{lem:pre-det facts}(4), together with (\ref{eq:canonical map is iso}), allows us to realize this map as precisely the pre-determinant applied to a loop.
\end{proof}
\end{lemma}
 
\subsection{Covering spaces}
Finally, let us say something about covering groups.
We use results from \cite[Chapter 10]{RotmanATBook}; one can also see \cite[Chapter 1.3]{HatcherAT}.
Covering spaces are common objects in algebraic topology, and are useful for computing homotopy groups of relatively nice spaces.
For example, one can use the universal covering space of the unit circle to show that its fundamental group can be identified with the integers.

A \emph{covering space} of a topological space $X$ consists of a topological space $\tilde{X}$ and a continuous surjective map $p: \tilde{X} \to X$ such that for each $x \in X$, there is an open neighbourhood $U$ of $x$ in $X$ such that $p^{-1}(U)$ is a union of disjoint open sets, each of which is mapped homeomorphically onto $U$ by $p$.
A covering space $\tilde{X}$ of $X$ is a \emph{universal covering space} if $\tilde{X}$ is simply connected; that is, it is path-connected with trivial fundamental group. 

  \begin{example}
  The real numbers $\bR$ is the universal covering space of the unit circle $\bT$ with covering map $p(t) = e^{2\pi i t}$. In fact, this generalizes to the following: $SU(n) \times \bR$   is a cover of $U(n)$.
  Indeed, we have covers
  \begin{equation}
  SU(n) \times \bR \to SU(n) \times \bT \to U(n)
  \end{equation}
  where the first map is $\id \times p$, where $p$ is the covering map $\bR \to \bT$ coming from the universal cover of $\bT$, and $SU(n) \times \bT \to U(n)$ is the map $(U,\lambda) \mapsto \lambda U$. 
  As $SU(n) \times \bR$ is simply connected, it follows that this is the universal cover of $U(n)$. 
  \end{example}  

  For a topological group $G$, a \emph{covering group} $\widetilde{G}$ of $G$ is a covering space of $G$ such that $\widetilde{G}$ is a topological group and the covering map $\rho: \widetilde{G} \to G$ is a continuous group homomorphism.
  A covering group is the \emph{universal covering group} if it is universal as a covering space. In fact, covering spaces of topological groups can be equipped with a \emph{some} multiplication which make them into a covering group \cite[Theorem 10.42]{RotmanATBook}.
  
  A topological space $X$ is said to be \emph{locally path-connected} if every neighbourhood of every point contains a path-connected neighbourhood and \emph{semilocally simply connected} if each point has a neighbourhood $N$ with the property that every loop in it can be contracted to a single point within $X$ (i.e., every loop in $N$ is null-homotopic in $X$). 
  
  \begin{theorem}
    A topological group $G$ which is path-connected, locally path-connected, and semilocally simply connected has a universal covering $\widetilde{G}$, unique up to homeomorphic isomorphism. Moreover, $\tilde{G}$ can be realized as
  \begin{equation}
  \widetilde{G} = PG/\sim_{nh}
  \end{equation}
  where $\sim_{nh}$ is the equivalence relation of null-homotopy. 
  \end{theorem}
  As mentioned in the introdction, the universal covering group of $U(n)$ is a semidriect product $SU(n) \rtimes \bR$.

  When a topological group $G$ has a universal covering group, we obtain a short exact sequence
  \begin{equation}\label{eq:covering group ses}
  \begin{tikzcd}
1 \arrow[r] & \pi_1(G) \arrow[r] & \widetilde{G} \arrow[r] & G \arrow[r] & 1
\end{tikzcd}
\end{equation}
where the first map is just the canonical embedding of homotopy classes into $\tilde{G}$, and the second map is evaluation at 1.

\section{The universal covering group of $U(\cM)$}\label{sec:universal cover of U(M)}

We have spoken about the universal covering group of $U(n)$ in the introduction, so let us move on to the unitary groups $U(\cM)$ where $\cM$ is not of type $I_n$.
In the case where $\cM$ is a properly infinite von Neumann algebra, the unitary group is path-connected and locally path-connected by the arguments below. For the last condition, one can appeal to \cite[Theorem 3.5]{Handelman78} to see that $\pi_1(U(\cM)) = 0$, i.e., that $U(\cM)$ is simply connected and therefore semilocally simply connected (as it is connected).
As such, the universal covering group of $U(\cM)$ exists, but it is trivially just $U(\cM)$ due to the fact that $\pi_1(U(\cM)) = 0$.

Let us first prove that the universal covering group of $U(\cM)$ exists whenever $\cM$ is a type $\II_1$ von Neumann algebra. 

\begin{lemma}
Let $\cM$ be a $\II_1$ von Neumann algebra. Then $U(\cM)$ is path-connected, locally path-connected, and semi-locally simply connected. In particular, the universal covering group of $U(\cM)$ exists and can be identified with
\begin{equation}
PU(\cM)/\sim_{nh}.
\end{equation}
\begin{proof}
Clearly $U(\cM)$ is path-connected as every unitary in a von Neumann algebra can be written as an exponential $u = e^{2\pi i a}$ for some self-adjoint $a \in A$ (by using functional calculus), and so the path $t \mapsto e^{2\pi i ta}$ yields a path from 1 to $u$.
The second condition amounts to showing that every open neighbourhood of an element contains a path-connected open neighbourhood.
This is trivial because if $\|u - v\| < 2$, then $u,v$ are path-connected (see, for example, \cite[Lemma 2.1.3(iii)]{RordamKBook}).
Thus if $N$ is a neighbourhood of $u \in U(\cM)$, let $0 < \ee < 1$ be such that $B_\ee(u) \subseteq N$.
Then any two points in $B_\ee(u)$ will be within distance 2 of each other.

The third condition says that every element $u \in U(\cM)$ has a neighbourhood $N$ such that every loop in $N$ is null-homotopic to a point in $U(\cM)$. To see this, let $u \in U(\cM)$ and consider $N = B_\frac{1}{2}(u)$, the open ball of radius $\frac{1}{2}$ around $u$. Take some $u_0 \in N$ and suppose that $\xi: [0,1] \to N$ is a loop with initial and terminal point $u_0$. 
We note that
\begin{equation}
\|\xi(t) - u_0\| \leq \|\xi(t) - u\| + \|u - u_0\| < \frac{1}{2} + \frac{1}{2} = 1.
\end{equation}
Define $\xi_0(t) = u_0^*\xi(t)$, which is a loop with initial and terminal point $1_\cM$, and note that
\begin{equation}
\begin{split}
\|\xi_0(t) - 1\| &= \|u_0^*\xi(t) - 1\| \\
&= \|\xi(t) - u_0\| < 1.
\end{split}
\end{equation}
We can now apply \ref{lem:pre-det facts}(2) to get that
\begin{equation}
\tD_\tau(\xi_0) = \frac{1}{2\pi i}\Tr_\cM\left(\log\xi_0(1)\right) = \frac{1}{2\pi i}\Tr_\cM(\log 1) = 0.
\end{equation}
In particular, as $\tD|_{\pi_1(U(\cM))}: \pi_1(U(\cM)) \to Z(\cM)_{sa}$ is an isomorphism, we have that $\xi_0$ is null-homotopic.
But then $\xi = u_0\xi_0$ is homotopic to the constant path $t \mapsto u_0$. 
\end{proof}
\end{lemma}

\begin{remark}
Suppose that $A$ is a unital C*-algebra and $\tD: PU^0(A) \to \Aff T(A)$ is the pre-determinant associated to the tracial map $\Tr_A: A_{sa} \to A_{sa}/A_0$ (which is just the quotient map) where $A_0$ is the subspace of all self-adjoint elements which vanish on every tracial state. $\Tr_A$ is often referred to as the universal trace.
If $\tD$ ``determines homotopy'' in the sense that whenever $\xi$ is a loop, we have that $\tD(\xi) = 0$ if and only if $\xi$ is null-homotopic, then the same argument above will yield that the universal covering group of $U^0(A)$ will exist. 

Thus whenever $A$ is $K_1$-injective, has trivial $K_1$ group, and satisfies the above property of the pre-determinant realizing homotopy, then $U(A)$ will have a universal cover.
For example, the unitary group of the Jiang-Su algebra $\cZ$ \cite{JiangSu99} or the unitary group of a UHF algebra will have a universal cover.
\end{remark}

Recall that a short exact sequence
\begin{equation}
\begin{tikzcd}
1 \arrow[r] & K \arrow[r, "\alpha"] & S \arrow[r, "\beta"] & U \arrow[r] & 1
\end{tikzcd}
\end{equation}
of groups splits if there is a group homomorphism $\sigma: U \to S$ such that $\beta \circ \sigma = \id_U$, and this happens if and only if $S = \alpha(K) \rtimes U'$, where $U' \leq S$ is some copy of $U$ (namely $\sigma(U)$) with $\beta|_{U'} \to U$ being an isomorphism. Moreover, this is in fact a direct product if and only if this copy $U'$ is normal in $S$. We will need the following well-known fact.

\begin{lemma}\label{lem:left split ses}
Let
\begin{equation}
\begin{tikzcd}
1 \arrow[r] & K \arrow[r, "\alpha"] & S \arrow[r, "\beta"] \arrow[l, "\gamma", bend left=49] & U \arrow[r] & 1
\end{tikzcd}
\end{equation}
be a short exact sequence of groups with $\gamma \circ \alpha = \id_K$. Then $S \simeq K \times U$.
\end{lemma}

\begin{theorem}
Let $\cM$ be a $\II_1$ von Neumann algebra and consider the universal covering group
\begin{equation}\label{eq:universal cover of U}
\begin{tikzcd}
1 \arrow[r] & \pi_1(U(\cM)) \arrow[r] & \widetilde{U(\cM)} \arrow[r] & U(\cM) \arrow[r] & 1.
\end{tikzcd}
\end{equation}
The sequence (\ref{eq:universal cover of U}) splits algebraically, giving
\begin{equation}
\widetilde{U(\cM)} \simeq \pi_1(U(\cM)) \times U(\cM).
\end{equation}
In particular,
\begin{equation}\label{eq:cg is Z times U}
\widetilde{U(\cM)} \simeq Z(\cM)_{sa} \times U(\cM).
\end{equation}
\begin{proof}
Let us denote the map $\pi_1(U(\cM)) \to \widetilde{U(\cM)}$ by $\iota$ and note that $\iota([\xi]) = [\xi]$, i.e., it just takes the homotopy class of $\xi$ to the class of $\xi$ under null-homotopy.
We will write $\ev_1: \widetilde{U(\cM)} \to U(\cM)$ to be evaluation at 1.  
By Lemma \ref{lem:fundamental group is center}, we can take the perspective that $K_0(\cM)$ is $\pi_1(U(\cM))$ with
\begin{equation}
\tD|_{\pi_1(U(\cM))}: \pi_1(U(\cM)) \to Z(\cM)_{sa}
\end{equation}
being an isomorphism.

Define $\gamma: \widetilde{U(\cM)} \to \pi_1(U(\cM))$ via the composition
\begin{equation}
\begin{tikzcd}
\widetilde{U(\cM))} \arrow[r, "\tD_0"] & Z(\cM)_{sa} \arrow[r, "(\tD|_{\pi_1(U(\cM))})^{-1}"] & \pi_1(U(\cM))
\end{tikzcd}
\end{equation}
where $\tD_0: \widetilde{U(\cM)} \to Z(\cM)_{sa}$ is the map given by $[\xi] \mapsto \tD(\xi)$, which is well-defined since $\tD$ preserves homotopy classes of paths with the same endpoints.
That is, given some class $[\xi]$ of a path $[0,1] \to U(\cM)$, up to null-homotopy, we evaluate the pre-determinant at this (homotopy class of a) path which gives us a self-adjoint element of the center. We then find a loop, up to homotopy, which gives that same self-adjoint element when we evaluate the pre-determinant. But then we realize that
\begin{equation}
\gamma \circ \iota([\xi]) = \gamma([\xi]) = [\xi] \text{ for } [\xi] \in \pi_1(U(\cM)),
\end{equation}
i.e., $\gamma \circ \iota = \id_{\pi_1(U(\cM))}$.
So we have the left split short exact sequence
\begin{equation}
\begin{tikzcd}
1 \arrow[r] & \pi_1(U(\cM)) \arrow[r, "\iota"] & \widetilde{U(\cM)} \arrow[r, "\ev_1"] \arrow[l, "\gamma", bend left=49] & U(\cM) \arrow[r] & 1,
\end{tikzcd}
\end{equation}
which gives that $\widetilde{U(\cM)} \simeq \pi_1(U(\cM)) \times U(\cM)$ by Lemma \ref{lem:left split ses}. Now (\ref{eq:cg is Z times U}) follows from our identification of $\pi_1(U(\cM))$ with $Z(\cM)_{sa}$.
\end{proof}
\end{theorem}

\begin{remark}
Even though we have a direct product as groups, it is clear that this is not a direct product as topological spaces. Indeed, a universal cover must be simply connected. This clearly fails to be true if $Z(\cM)_{sa}$ and $U(\cM)$ were equipped with their (relative) norm topologies. 
Moreover, $\pi_1(U(\cM)) \subseteq \widetilde{U(\cM)}$ is considered as a discrete set.
\end{remark}

Recall that a group $G$ is perfect if $G = DG$, i.e., every element is a finite product of commutators.
With the above, we conclude that the universal covering group of $U(\cM)$ is not perfect.

\begin{cor}
The universal covering group of $U(\cM)$ is not perfect whenever $\cM$ is a $\II_1$ von Neumann algebra.
\begin{proof}
If $G = A \times H$ where $H$ is some group and $A$ is abelian, we have that $DG = DA \times DH = \{0\} \times DH$. So if $|A| > 1$, then $G$ can never be perfect. 
\end{proof}
\end{cor}

  \bibliographystyle{amsalpha}
  \bibliography{biblio}

\end{document}